\documentclass[letterpaper, 10 pt, conference]{ieeeconf}

\usepackage{amssymb,amsmath,color}
\usepackage{graphicx}
\usepackage{caption}
\usepackage{enumerate}
\IEEEoverridecommandlockouts                              % This command is only
\newcommand*{\QEDB}{\hfill\ensuremath{\square}}%

\newtheorem{assumption}{Assumption}
\newtheorem{definition}{Definition}

\newtheorem{theorem}{Theorem}

\newtheorem{proposition}{Proposition}
\newtheorem{remark}{Remark}

\newtheorem{example}{Example}

\begin{document}

\title{\bf Fixed-Time ISS for Interconnected Discrete-Time Systems, Part II:\\ Finite-Horizon Lyapunov Functions}

\author{Sumadhu Rubaiyat \thanks{S. Rubaiyat and J. I. Poveda are with the ECE Department, UC, San Diego. Email:\{poveda@ucsd.edu.\}. Research supported in part by NSF grant CAREER
2305756.} and Jorge I. Poveda.}

\maketitle
%%%%%%%%%%%%%
\begin{abstract}
We study fixed-time input-to-state stability (FxT-ISS) for discrete-time systems subject to exogenous inputs using finite-horizon Lyapunov functions. In the proposed Lyapunov-based characterization, the effect of the initial condition is eliminated after a uniformly bounded number of iterations. This characterization encompasses existing extinction-based Lyapunov conditions while allowing more general finite-memory mechanisms. Building on this framework, we establish a finite-horizon small-gain theorem for interconnected discrete-time systems, providing conditions under which FxT-ISS is preserved under feedback interconnections. The proposed framework is illustrated through examples arising in control, distributed algorithms, and feedback optimization.
\end{abstract}
%%%%%%%%%%%%%%%%%%%%%%%%%%%%%%%%%%%%%%%%%%
\section{INTRODUCTION}
Traditionally, discrete-time control, optimization, and decision-making algorithms have been analyzed using tools that certify \emph{asymptotic stability} \cite{jiang2001input,Kellet_inclusions1,geiselhart2017equivalent,Haddad_Nonlinear}. This property combines stability, which prevents the iterates from undergoing arbitrarily large excursions from their nominal operating point, with asymptotic attractivity, which guarantees convergence to that point as time tends to infinity. However, increasing computational and performance demands have motivated the development of algorithms that converge not merely asymptotically, but exactly within a uniformly bounded time. This behavior is captured by the notion of \emph{fixed-time stability}. In the continuous-time setting, since its introduction in \cite{andrieu2008homogeneous}, fixed-time stability has received significant attention, leading to new algorithmic architectures for control \cite{lopez2019conditions,polyakov2011nonlinear}, estimation \cite{rios2017time}, learning \cite{poveda2022fixed,poveda2021nonsmooth}, and optimization \cite{garg2022fixed}. By comparison, fixed-time stability remains considerably less developed in the \emph{discrete-time} setting.

Finite-step Lyapunov functions have been developed for discrete-time systems as relaxations of standard ``one-step'' Lyapunov conditions \cite{noroozi2020control}. For example, the finite-step ISS framework presented in \cite{geiselhart2016relaxed} replaces one-step dissipation by a contraction condition imposed after a fixed number of iterations, while retaining a dependence on the initial Lyapunov value. This framework has subsequently enabled relaxed small-gain results for discrete-time interconnections. In parallel, fixed-time stability of discrete-time systems has been studied through Lyapunov difference inequalities specifically designed to guarantee finite convergence with a settling-time bound independent of the initial condition \cite{lee2023fixed,tatari2023deterministic}.  Finite-time ISS and associated small-gain conditions for interconnected discrete-time systems have also been investigated in \cite{mo2022small}, where finite-time ISS requires the state to become bounded by a function of the input after a finite time. More recently, in
\cite{mazenc2024fixed}, the authors introduced a fixed-time ISS notion for discrete-time systems characterized by a finite-memory property: after a uniform transient, the state is bounded exclusively in terms of disturbances acting over a fixed-length recent time window. They developed trajectory-based conditions and sampled-data controllers achieving this property for time-varying linear systems, including extensions to output feedback and input delays. Finite-time results for discrete-time systems were also studied in \cite{mo2024lyapunov} using a small-gain theorem.

Motivated by the above background, this paper makes two main theoretical contributions. First, we introduce a novel finite-horizon Lyapunov framework for certifying FxT-ISS in discrete-time systems, in which the Lyapunov estimate completely loses its dependence on the initial state after a uniformly bounded number of iterations. Unlike existing finite-step ISS Lyapunov functions, which enforce contraction over a finite number of steps, the proposed condition requires the initial-state contribution to disappear entirely. Also, unlike existing fixed-time conditions based on specific one-step dissipation laws, it allows fixed-time convergence to emerge from the multi-step structure of the dynamics. We establish the corresponding FxT-ISS theorem and show that the framework strictly encompasses the extinction-based Lyapunov conditions studied in \cite{lee2023fixed,tatari2023deterministic}, as well as the companion paper \cite{ACC26A}. Second, inspired by recent continuous-time small-gain results for FxT stability \cite{tang2026lyapunov,tang2026smallgain}, we develop a finite-horizon small-gain theorem for interconnected discrete-time systems, where the influence of the initial subsystem states is eliminated after finitely many compositions of the finite-horizon gains, rather than merely contracted asymptotically. The proposed framework thereby captures a range of distinct finite-horizon mechanisms, which we illustrate through deadbeat dynamics, finite-information propagation in min-consensus, and feedback optimization. The results are also contrasted with the extinction-based Lyapunov approach to FxT-ISS studied in the companion paper \cite{ACC26A}, showing that, when the input is set to zero, continuous nilpotent update maps admit finite-horizon certificates, while those with unbounded images admit no extinction-based certificate relative to the origin. 

The rest of this paper is organized as follows. Section \ref{sec_notation} presents the preliminaries. Sections \ref{mainresults}  and IV present the main results, Section \ref{secproofs} presents the proof, and Section \ref{sec:conclusions} ends with the conclusions.

%%%%%%%%%%%%%%%%%%%%%%%%%%%%%%%%%%%%%%%
\section{PRELIMINARIES}
\label{sec_notation}
\subsection{Notation}
Given a nonempty closed set $\mathcal{A} \subset \mathbb{R}^n$ and a vector $z \in \mathbb{R}^n$, we define $|z|_\mathcal{A} := \min_{y \in \mathcal{A}}|z-y|$, and we use $|\cdot|$ to denote the standard Euclidean norm.  Throughout the paper, for ease of notation we write $(u,v)$ for $(u^{T},v^{T})^{T}$.  We let $\mathcal{PD}$ be the class of continuous functions $\rho:\mathbb{R}_{\geq0}\to\mathbb{R}_{\geq0}$ that are positive definite. A continuous function $\alpha:\mathbb{R}_{\geq0}\to\mathbb{R}$ is said to be of class $\mathcal{K}$ if it is zero at zero and strictly increasing, and it is of class $\mathcal{K}_{\infty}$ if, additionally, it is unbounded. A function $\beta:\mathbb{R}_{\geq 0}\times\mathbb{Z}_{\geq 0}
\to\mathbb{R}_{\geq 0}$ belongs to $\mathcal{KL}_{\mathrm{FxT}}$ if there exists
$M\in\mathbb{Z}_{\geq 1}$ such that: (a) for every $k\in\mathbb{Z}_{\geq 0}$, $\beta(\cdot,k)$ is continuous, nondecreasing, and $\beta(0,k)=0$; (b) for every $k<M$, $\beta(\cdot,k)$ belongs to class $\mathcal{K}$; (c) for every $r\geq0$, $\beta(r,\cdot)$ is  nonincreasing; (d) $\beta(r,k)=0$ for every $r\geq0$ and every $k\geq M$. For any $r\in\mathbb{R} $, we use $[r]_+:=\max\{r,0\}$, and for any map $h:\mathbb R_{\geq0}\to\mathbb R_{\geq0}$,
we use $h^k$ to denote its $k$-fold composition, defined
recursively by $h^0:=\operatorname{id}$, $h^{k+1}:=h\circ h^k$.
\section{FINITE-HORIZON LYAPUNOV-BASED CONDITIONS FOR FXT-ISS}
\label{mainresults}
In this paper, we consider time-invariant discrete-time dynamical systems with inputs, given by
\begin{equation}\label{eq:dynamics}
x^+=g(x,u),~~~~x_0:=x(0)\in\mathbb{R}^n,
\end{equation}
where $x\in\mathbb R^n$, $u\in\mathbb R^m$,
and $g:\mathbb R^n\times\mathbb R^m\to\mathbb R^n$ is continuous. Let $\mathcal A\subset\mathbb R^n$ be
a nonempty closed set that is forward invariant for the nominal system $x^+=g(x,0)$. For this set, we study the following property, which is a natural extension to the discrete-time setting of similar FxT-ISS notions studied for continuous-time systems \cite{lopez2019conditions}. 
%
%\begin{assumption}\label{assumption1}
%The function $g(\cdot,\cdot)$ is continuous, and the set $\mathcal{A}$ is nonempty and closed. \hfill $\QEDB$
%\end{assumption}
%

%In this paper, we will study the following notion of stability under inputs/disturbances:

\vspace{0.1cm}
\begin{definition}\label{def:fxISS}
System \eqref{eq:dynamics} is said to render the set $\mathcal{A}$ fixed-time input-to-state stable (ISS) if there exists a class $\mathcal{K}\mathcal{L}_{\text{FxT}}$ function $\beta$ and a class $\mathcal{K}$ function $\gamma(\cdot)$, such that for all $x_0\in\mathbb{R}^n$ and all inputs $u\in\ell^{\infty}$, all solutions of \eqref{eq:dynamics} satisfy the following inequality
\begin{equation}\label{eq:general_fxtiss_estimate}
    |x(j)|_{\mathcal A}\leq\max\left\{\beta(|x(0)|_{\mathcal A},j),\gamma\left(|u|_{[0,j-1],\infty}\right)\right\},
\end{equation}
for all $j\geq0$, where $|u|_{[0,j-1],\infty}=0$ for $j=0$.  \hfill $\square$
\end{definition}
\subsection{Sufficient Finite-Horizon Lyapunov-based Conditions}
We will study the FxT-ISS property using the following class of Lyapunov functions.
\begin{definition}
A continuous function
$V:\mathbb R^n\to\mathbb R_{\geq0}$ is said to be a \emph{Finite-Horizon FxT-ISS Lyapunov function} for system \eqref{eq:dynamics} if there exist functions $\alpha_1,\alpha_2,\sigma,\gamma
\in\mathcal K_\infty$, a constant $\rho\in[0,1)$, and an integer
$M\geq1$ such that the following
conditions hold:
\begin{enumerate}[(a)]
\item  For every $x\in\mathbb R^n$,
\begin{equation}
    \alpha_1(|x|_{\mathcal A})
    \leq V(x)
    \leq \alpha_2(|x|_{\mathcal A}).
    \label{eq:fh_comparison}
\end{equation}
\item For every $(x,u)\in\mathbb R^n
\times\mathbb R^m$,
\begin{equation}
    V(x^+)-V(x)
    \leq -(1-\rho) V(x)+\sigma(\|u\|).
    \label{eq:fh_dissipation}
\end{equation}
\item Let $G_j$ be the $j$-step transition map defined recursively with $G_0(x) := x$ as 
\[
G_{j+1}(x,u_0,\ldots,u_j)
:=
g\!\left(G_j(x,u_0,\ldots,u_{j-1}),u_j\right).
\]
Then, the following bound holds
\begin{equation}\label{finitehorizonbound}
V\!\left(G_M(x,u_0,\ldots,u_{M-1})\right)
\leq
\gamma\!\left(
\max_{0\leq j\leq M-1}\|u_j\|
\right),
\end{equation}
for every $x\in\mathbb{R}^n$ and every
$u_0,\ldots,u_{M-1}\in\mathbb{R}^m$, 
\end{enumerate}
\end{definition}
\vspace{0.1cm}
\begin{proposition}\label{equivalentprop}
Consider the discrete-time system \eqref{eq:dynamics}, where $g:\mathbb R^n\times\mathbb R^m\to\mathbb R^n$ is continuous,
and let $\mathcal A\subset\mathbb R^n$ be compact.  Then a continuous function $V:\mathbb R^n\to\mathbb R_{\geq0}$ is a finite-horizon FxT-ISS Lyapunov function if and only if
there exist $\alpha_1,\alpha_2,\chi,\gamma\in\mathcal K_\infty$,
a constant $\rho\in[0,1)$, and an integer $M\geq1$ such that \eqref{eq:fh_comparison} and \eqref{finitehorizonbound} hold, as well as
\begin{equation}\label{implicationcondition}
V(x)\geq\chi(\|u\|)
\quad\Longrightarrow\quad
V(g(x,u))\leq\rho V(x).
\end{equation}
for all $(x,u)\in\mathbb R^n\times\mathbb R^m$.
\end{proposition}

\vspace{0.1cm}
The following theorem is the first main result of this paper. All proofs are presented in Section \ref{secproofs}.

\vspace{0.1cm}
\begin{theorem}\label{theoremLyapu}
Suppose that $V$ is a finite-horizon FxT-ISS Lyapunov function. Then $\mathcal A$ is FxT-ISS. \hfill $\QEDB$
%where $\|u\|_{[0,k-1],\infty}=0$ at $k=0$, and
%\begin{equation}
%    \beta(s,k)=
%    \begin{cases}
%    \displaystyle
%    \alpha_1^{-1}
%    \left(
%        2\rho^k\alpha_2(s)
%    \right),
%    &0\leq k<M,\\[6pt]
%    0,
%    &k\geq M,
%    \end{cases}
%    \label{eq:fh_beta}
%\end{equation}
%and
%\begin{equation}
%    \chi(s)=
%    \alpha_1^{-1}
%    \left(
%        2\max\left\{
%            \frac{\sigma(s)}{1-\rho},
%            \gamma(s)
%        \right\}
%    \right).
%    \label{eq:fh_gain}
%\end{equation}
\end{theorem}
\vspace{0.1cm}
\begin{remark}\label{rem:finite_phase}
Theorem~\ref{theoremLyapu} also applies to systems of the form:
\[
z^+=g(z,q,u),~~q^+=\tau(q),~~(z,q)\in \mathbb R^d\times\mathcal Q,
\]
where $\mathcal Q$ is a nonempty finite set endowed with the discrete topology and $\tau:\mathcal Q\to\mathcal Q$ is any
map. The Lyapunov inequalities are then required to hold
along the resulting augmented dynamics, uniformly over all initial indices $q(0)\in\mathcal Q$, with common gains, a common horizon, and $\mathcal{A}=\mathcal A_0\times\mathcal Q$. \QEDB
\end{remark}

\subsection{Connections with Extinction-based FxT-ISS Lyapunov Functions:}
In \cite{ACC26A}, the FxT-ISS property was studied using the so-called \emph{Extinction-based FxT-ISS Lyapunov functions}. These functions leverage a class of dissipation maps $\psi\in\mathcal{PD}$ that satisfy the property $h_{\psi}^M(s)=0$ for some $M\in\mathbb{Z}_{\geq1}$ and all $s\geq0$, where $h_{\psi}(s)=\max_{0\leq r\leq s}[r-\psi(r)]_+$. Functions of this form are said to be of class $\mathcal{E}_{\text{FxT}}^M$. In particular, the following definition corresponds to \cite[Def. 3]{ACC26A}.

\vspace{0.1cm}
\begin{definition}\cite[Def. 3]{ACC26A} \label{def:general_fxtiss_lyap}
A continuous function $V:\mathbb R^n\to\mathbb R_{\geq0}$ is called an \emph{Extinction-based FxT-ISS
Lyapunov function} with respect to $\mathcal A\subset\mathbb{R}^n$ if there exist $\alpha_1,\alpha_2\in\mathcal K_\infty$, $\rho\in\mathcal K_\infty$, and $\psi\in\mathcal{E}^M_{\text{FxT}}$, such that:
\begin{subequations}\label{eq:suffLyap1}
\begin{align}   
    \alpha_1(|x|_{\mathcal A})&\leq V(x)\leq \alpha_2(|x|_{\mathcal A}),
    \qquad \forall x\in\mathbb R^n, \label{eq:general_lyap_bounds}\\
    V(x^+)-V(x)&\leq-\min\{\psi(V(x)),V(x)\}+\rho(|u|)\label{eq:additive_dissipation}
\end{align}
\end{subequations}
for every $(x,u)\in\mathbb R^n\times\mathbb R^m$. \hfill $\square$
\end{definition}

As shown in \cite[Thm. 1]{ACC26A}, Extinction-based FxT-ISS
Lyapunov functions certify the FxT-ISS property by evaluating the decrease of the Lyapunov function after every jump of \eqref{eq:dynamics}, and without the need to evaluate the map $G_M$ used in \eqref{finitehorizonbound}. In fact, the following proposition shows that any Extinction-based FxT-ISS
Lyapunov is also a Finite-Horizon FxT-ISS Lyapunov function. The converse, however, is in general not true.
\begin{proposition}\label{thm:extinction_implies_fh}
Let $\mathcal{A}\subset\mathbb{R}^n$ be a nonempty closed set, and suppose that $V:\mathbb{R}^n\to\mathbb{R}_{\geq 0}$ is an
extinction-based FxT-ISS Lyapunov function for system \eqref{eq:dynamics} relative to $\mathcal{A}$. Then $V$ is also a finite-horizon FxT-ISS Lyapunov function for system \eqref{eq:dynamics} relative to $\mathcal{A}$. \hfill $\QEDB$
\end{proposition}

\begin{example} The class of finite-horizon FxT-ISS Lyapunov functions is strictly broader than the class of extinction-based FxT-SS Lyapunov functions. To see this, we can focus on the case when $u=0$ and consider the discrete-time linear system $x^+=Ax$ with $A=[0,a;0,0]$ (using ``Matlab" notation) with $a\in(0,1)$, with $\mathcal{A}=\{0\}$. Consider the Lyapunov function $V(x)=\|x\|$, which clearly satisfies \eqref{eq:fh_comparison}. Since $\|A\|=a$, it follows that $ V(x^{+})\leq aV(x)$, and therefore \eqref{eq:fh_dissipation} holds with zero input gain. Moreover, $A^2=0$, and hence every solution satisfies $ x(j+2)=0$ for every $j\in\mathbb{Z}_{\geq0}$. Therefore, $V(x(j+2))-V(x(j)) =-V(x(j))$, which implies that \eqref{finitehorizonbound} holds with $M=2$. Therefore, $V$ is a finite-horizon FxT (ISS) Lyapunov function. Next, we show that $V$ cannot be an extinction-based FxT Lyapunov
function. Suppose, by contradiction, that there exist
$N\in\mathbb{Z}_{>0}$ and
$\psi\in\mathcal{E}_{\mathrm{FxT}}^{N}$ such that \eqref{eq:additive_dissipation} holds with $\rho=0$. For any $s>0$, let $x=(0,s)$. Then, $V(x)=s$, $x^{+}=[as;0]$, and  $V(x^{+})=as$. Substituting these expressions into \eqref{eq:additive_dissipation} gives $-(1-a)s \leq-\min\{\psi(s),s\}$, or, equivalently, $\min\{\psi(s),s\}\leq(1-a)s$. Since $a\in(0,1)$, we have $(1-a)s<s$ for every $s>0$. Thus, \eqref{eq:additive_dissipation} necessarily implies $\psi(s)\leq(1-a)s$, for all $s>0$. Let $g(s):=[s-\psi(s)]_{+}$. It follows that $g(s)\geq as$ for all $s>0$. Therefore, its monotone envelope $h(s):=\max_{0\leq r\leq s}g(r)$ satisfies $h(s)\geq as$ for all $s>0$. In particular, $h$ is unbounded. This contradicts
$\psi\in\mathcal{E}_{\mathrm{FxT}}^{N}$, since the fixed-time
extinction property requires $h^{N}(s)=0$, for all $s\geq0$, which, in particular, implies that $h$ is globally bounded. Hence, there is no $\psi\in\mathcal{E}_{\mathrm{FxT}}^{N}$, for any finite $N$, such that \eqref{eq:additive_dissipation} holds. \hfill $\square$
\end{example}

The following proposition complements the previous results for the case when $u=0$, and provides further justification for the study of finite-horizon FxT (ISS) Lyapunov functions for discrete-time systems.
\begin{proposition}\label{propnecessarycondition}
 Let $\mathcal A\subset\mathbb{R}^n$ be a nonempty closed set, and consider the system $x^+=g(x)$, where $g:\mathbb{R}^n\to\mathbb{R}^n$ is continuous and $g(\mathcal{A})\subset\mathcal{A}$. Then, the following holds:
\begin{enumerate}[(a)]
\item If the system admits an extinction-based FxT-ISS Lyapunov function relative to $\mathcal A$, then
\[
\sup_{x\in\mathbb{R}^n}|g(x)|_{\mathcal A}<\infty.
\]
\item Let $\mathcal A=\{0\}$ and suppose there exists $M\in\mathbb{Z}_{\geq1}$ such that $g^M(x)=0$ for all $x\in\mathbb{R}^n$. Then the system admits a finite-horizon FxT-ISS Lyapunov function.  \hfill $\square$
\end{enumerate}
\end{proposition}
\begin{remark}
A direct consequence of Propostion \ref{propnecessarycondition} is that if $\sup_{x\in\mathbb{R}^n}|g(x)|_{\mathcal A}=\infty$, then the dynamics $x^+=g(x)$ do not admit an extinction-based FxT-ISS Lyapunov function relative to $\mathcal{A}$. \QEDB
\end{remark}
\subsection{Examples in Control, Multi-Agent Coordination, and Optimization}
We illustrate Theorem \ref{theoremLyapu} via different examples that admit a finite-horizon FxT-ISS Lyapunov function.
\begin{example}[Deadbeat control under additive disturbances]
Consider the discrete-time linear system
\begin{equation}
\label{eq:deadbeat_system}
x^+=Ax+Bu,
\end{equation}
where $x\in\mathbb R^n$, $u\in\mathbb R^m$, and suppose that $A$ is nilpotent, i.e., there exists $M\in\mathbb Z_{\geq1}$ such that $A^M=0$. Such dynamics arise in deadbeat control. Indeed, for a controllable system
\[
z^+=A_0z+B_0v,
\]
a state-feedback law $v=Kz+u$ can be selected so that
$A:=A_0+B_0K$ has all its eigenvalues at the origin and is therefore
nilpotent. Let $\mathcal A=\{0\}$ and fix any $\rho\in(0,1)$. Define $V(x):=\sum_{\ell=0}^{M-1}\rho^{-\ell}\|A^\ell x\|$. We show that $V$ is a finite-horizon Lyapunov function for FxT-ISS. First, since the term corresponding to $\ell=0$ is $\|x\|$, then $V(x)\geq\|x\|$. Moreover, $V(x)\leq\left(\sum_{\ell=0}^{M-1}\rho^{-\ell}\|A^\ell\|\right)\|x\|$. Defining $c_A:=\sum_{\ell=0}^{M-1}
\rho^{-\ell}\|A^\ell\|$, the comparison bounds \eqref{eq:fh_comparison} hold with $\alpha_1(s)=s$ and $\alpha_2(s)=c_As$. Next, using the fact that $A$ is nilpotent
\begin{align*}
V(Ax)&=
\sum_{\ell=0}^{M-1}
\rho^{-\ell}\|A^{\ell+1}x\|=
\rho
\sum_{\ell=1}^{M}
\rho^{-\ell}\|A^\ell x\|
\\
&=
\rho
\sum_{\ell=1}^{M-1}
\rho^{-\ell}\|A^\ell x\|\leq
\rho V(x).
\end{align*}
Using \eqref{eq:deadbeat_system} and the triangle
inequality gives
\begin{align*}
V(x^+)
&\leq
V(Ax)+V(Bu)\leq
\rho V(x)
+
\sum_{\ell=0}^{M-1}
\rho^{-\ell}\|A^\ell Bu\|
\\
&\leq
\rho V(x)
+
\left(
\sum_{\ell=0}^{M-1}
\rho^{-\ell}\|A^\ell B\|
\right)\|u\|.
\end{align*}
Hence, defining $c_B:=\sum_{\ell=0}^{M-1}\rho^{-\ell}\|A^\ell B\|$, we obtain that \eqref{eq:fh_dissipation} holds:
\begin{equation}
\label{eq:deadbeat_iss_dissipation}
V(x^+)
\leq
\rho V(x)+\sigma(\|u\|),
\qquad
\sigma(r):=(1+c_B)r.
\end{equation}
Finally, to prove \eqref{finitehorizonbound}, we iterate \eqref{eq:deadbeat_system} for $M$ steps, obtaining
\begin{equation}
x(M)=A^Mx(0)+\sum_{j=0}^{M-1}A^{M-1-j}Bu(j).
\end{equation}
Since $A^M=0$, it follows that $x(M)
= \sum_{j=0}^{M-1}A^{M-1-j}Bu(j).$
%
%\begin{equation}
%\label{eq:deadbeat_M_state}
%x(M)
%=
%\sum_{j=0}^{M-1}
%A^{M-1-j}Bu(j).
%\end{equation}
Consequently,
\begin{align*}
V(x(M))
&\leq
\sum_{j=0}^{M-1}
V\!\left(A^{M-1-j}Bu(j)\right)
\\
&\leq
\left(
\sum_{j=0}^{M-1}
\sum_{\ell=0}^{M-1}
\rho^{-\ell}
\|A^{\ell+M-1-j}B\|
\right)
\|u\|_{[0,M-1],\infty}.
\end{align*}
Using $c_M:=
\sum_{j=0}^{M-1}
\sum_{\ell=0}^{M-1}
\rho^{-\ell}
\|A^{\ell+M-1-j}B\|$, we obtain $V(x(M))\leq c_M\|u\|_{[0,M-1],\infty}$. Since the same argument applies from any initial time, we obtain:
\begin{equation}
V\!\left(
G_M(x,u_0,\ldots,u_{M-1})
\right)
\leq
\gamma\!\left(
\max_{0\leq j\leq M-1}\|u_j\|
\right),
\end{equation}
with $\gamma(r):=(1+c_M)r$.
These choices ensure $\sigma,\gamma\in\mathcal K_\infty$
even when $B=0$.
Hence, $V$ is a finite-horizon Lyapunov function.
By Theorem~\ref{theoremLyapu}, the set $\mathcal A$
is FxT-ISS. \QEDB
\end{example}

\vspace{0.1cm}
The following example illustrates the application of Theorem \ref{theoremLyapu} in the well-known min-consensus algorithm.

\vspace{0.1cm}
\begin{example}[Min-Consensus under Perturbations]
Consider a strongly connected directed graph
$\mathcal G=(\mathcal V,\mathcal E)$, with
$\mathcal V=\{1,\ldots,N\}$ and $N\geq2$.
Let $D:=\max_{i,j\in\mathcal V}d(i,j)\geq1$
denote its directed diameter, where $d(i,j)$ is
the length of the shortest directed path from
node $i$ to node $j$. Each agent $i$ has a scalar state $x_i\in\mathbb R$ and evolves according to the dynamics:
\begin{equation}
\label{eq:min_consensus}
x_i^+
=
\min_{j\in\mathcal N_i\cup\{i\}} x_j+u_i,
\qquad i\in\mathcal V,
\end{equation}
where $u_i\in\mathbb R$ is an additive disturbance and $\mathcal{N}_i$ denotes the set of neighbors of agent \(i\). Define the nominal min-consensus map $T:\mathbb R^N\to\mathbb R^N$ componentwise as $T_i(x):=\min_{j\in\mathcal N_i\cup\{i\}}x_j$. Then \eqref{eq:min_consensus} can be written compactly as:
\begin{equation}\label{overallsystem}
x^+=T(x)+u.
\end{equation}
Let $\mathcal A:=\{c\mathbf 1:c\in\mathbb R\}$ be the consensus set, and define the disagreement measure $W(x):=\max_i x_i-\min_i x_i$. Since $\mathcal G$ is strongly connected and has directed diameter $D$, every node receives information originating from every other node after at most $D$ iterations. Hence,
\begin{equation}
\label{eq:min_consensus_D}
T^D(x)
=
\left(\min_i x_i\right)\mathbf 1,
\end{equation}
and therefore $W(T^D(x))=0$, for all $x\in\mathbb R^N$.  Fix any $\rho\in(0,1)$ and let
\begin{equation}
\label{eq:min_consensus_V}
V(x)
:=
\sum_{\ell=0}^{D-1}
\rho^{-\ell}W(T^\ell(x)).
\end{equation}
We show that $V$ is a finite-horizon Lyapunov function for system \eqref{overallsystem}, with horizon $M=D$. To see this, first note that
$W(T(x))\leq W(x)$ and therefore $W(T^\ell(x))\leq W(x)$ for all
$\ell\geq0$. Therefore, $W(x)\leq V(x)\leq c_DW(x),$
%
%\begin{equation}
%\label{eq:min_V_W_bounds}
%W(x)
%\leq V(x)
%\leq c_DW(x),
%\end{equation}
where $c_D:=\sum_{\ell=0}^{D-1}\rho^{-\ell}$. Since $|x|_{\mathcal A}=\min_{c\in\mathbb R}\|x-c\mathbf 1\|_2$, we have $W(x)\leq\sqrt{2}\,|x|_{\mathcal A}$. Moreover, if $\bar x:=\frac{1}{N}\mathbf 1^\top x$, then $|x|_{\mathcal A}=
\|x-\bar x\mathbf 1\|_2
\leq
\sqrt{N}\,W(x)$, so that $\frac{1}{\sqrt N}|x|_{\mathcal A}
\leq
W(x)
\leq
\sqrt{2}|x|_{\mathcal A}$. Combining the above inequalities gives $\frac{1}{\sqrt N}|x|_{\mathcal A}
\leq
V(x)
\leq
\sqrt{2}c_D|x|_{\mathcal A}$. Thus, $V$ satisfies \eqref{eq:fh_comparison} with $\alpha_1(s)=\frac{s}{\sqrt N}$ and $\alpha_2(s)=\sqrt{2}c_Ds$. Next, using \eqref{eq:min_consensus_D}, we get
\begin{align*}
V(T(x))&=
\sum_{\ell=0}^{D-1}
\rho^{-\ell}W(T^{\ell+1}(x))=
\rho
\sum_{\ell=1}^{D}
\rho^{-\ell}W(T^\ell(x))
\\
&=
\rho
\sum_{\ell=1}^{D-1}
\rho^{-\ell}W(T^\ell(x))\leq
\rho V(x).
\end{align*}
Since the min-consensus map is nonexpansive in the infinity norm, we have $\|T(x)-T(y)\|_\infty\leq\|x-y\|_\infty$. Also, we have that $|W(x)-W(y)|\leq2\|x-y\|_\infty$. Using the nonexpansiveness of $T$, it follows that, for every $\ell\geq0$,
\begin{align*}
W(T^\ell(T(x)+u))
&\leq
W(T^\ell(T(x)))
+
2\|T^\ell(T(x)+u)\\
&~~~~-T^\ell(T(x))\|_\infty
\\
&\leq
W(T^{\ell+1}(x))+2\|u\|_\infty.
\end{align*}
Consequently,
\begin{align*}
V(T(x)+u)
&=
\sum_{\ell=0}^{D-1}
\rho^{-\ell}W(T^\ell(T(x)+u))
\\
&\leq
\sum_{\ell=0}^{D-1}
\rho^{-\ell}W(T^{\ell+1}(x))
+
2\sum_{\ell=0}^{D-1}\rho^{-\ell}\|u\|_\infty
\\
&=
V(T(x))+2c_D\|u\|_\infty\leq
\rho V(x)+2c_D\|u\|_\infty.
\end{align*}
Thus, inequality \eqref{eq:fh_dissipation} holds with $\sigma(r):=2c_Dr$. Finally, we establish condition \eqref{finitehorizonbound}. Let $\bar x(0)=x(0)$, $x$ be a solution to \eqref{overallsystem}, and let $\bar x$ be a solution to \eqref{overallsystem} with $u=0$. By nonexpansiveness of $T$, we get
\begin{align*}
\|x(j+1)-\bar x(j+1)\|_\infty
&\leq
\|T(x(j))-T(\bar x(j))\|_\infty+\|u(j)\|_\infty
\\
&\leq
\|x(j)-\bar x(j)\|_\infty+\|u(j)\|_\infty.
\end{align*}
Iterating this inequality gives
\[
\|x(D)-\bar x(D)\|_\infty
\leq
\sum_{j=0}^{D-1}\|u(j)\|_\infty
\leq
D\|u\|_{[0,D-1],\infty}.
\]
By \eqref{eq:min_consensus_D}, we have $\bar x(D)
=\left(\min_i x_i(0)\right)\mathbf 1\in\mathcal A$. Since $W(\bar x(D))=0$, we get
\begin{align*}
W(x(D))&=|W(x(D))-W(\bar x(D))|\leq
2\|x(D)-\bar x(D)\|_\infty\\
&\leq
2D\|u\|_{[0,D-1],\infty}.
\end{align*}
Using the bounds on $W$, we get $V(x(D))\leq c_DW(x(D)) \leq 2Dc_D\|u\|_{[0,D-1],\infty}$. Since the same argument applies from any initial time $k$, we obtain
\[
V\!\left(
G_D(x,u_0,\ldots,u_{D-1})
\right)
\leq
\gamma\!\left(
\max_{0\leq j\leq D-1}\|u_j\|_\infty
\right),
\]
where $\gamma(r):=2Dc_Dr$. Hence, $V$ satisfies all the conditions of a finite-horizon FxT-ISS Lyapunov function. By Theorem \ref{theoremLyapu}, the set $\mathcal{A}$ is FxT-ISS. \QEDB 
\end{example}
\section{A FINITE-HORIZON LYAPUNOV-BASED SMALL-GAIN THEOREM}
\label{sec:smallgain}
In this section, we build on the finite-horizon Lyapunov-based
FxT-ISS characterization developed in the previous section to establish a fixed-time small-gain theorem for interconnected discrete-time systems of the form
\begin{figure}[t]
    \centering
    \includegraphics[width=0.9\columnwidth]{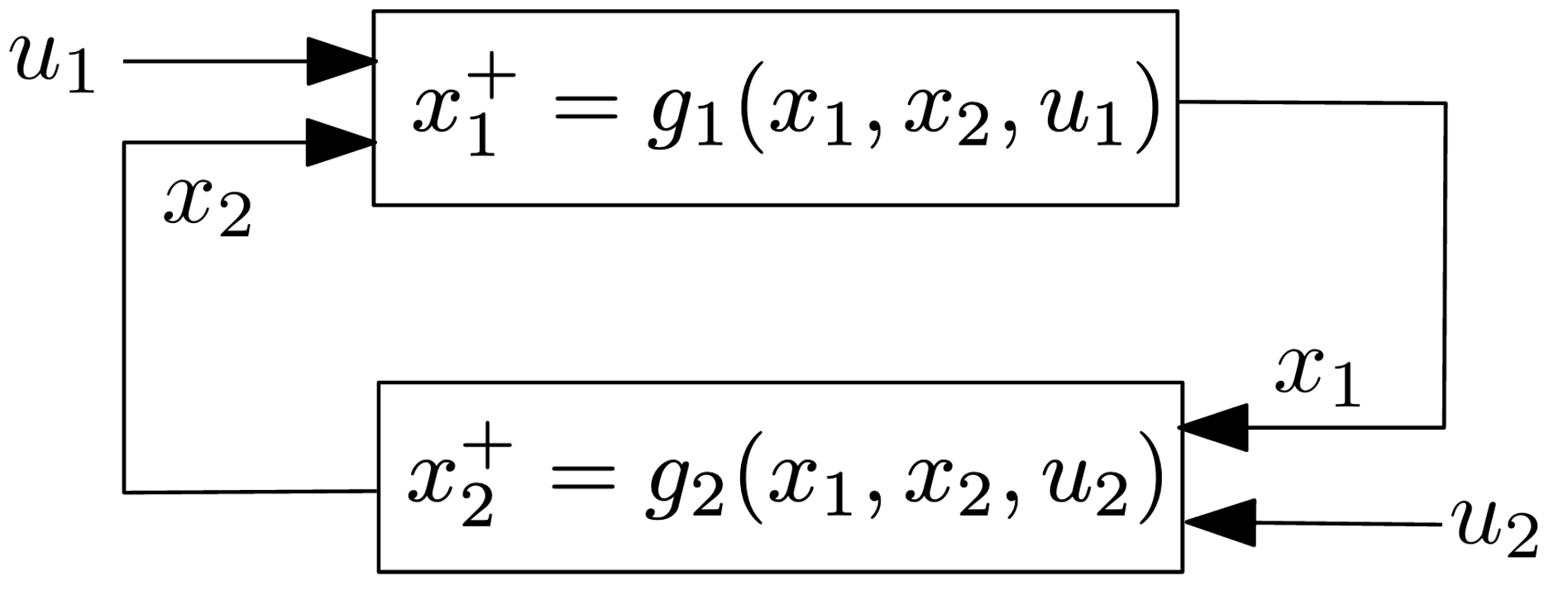}
    \caption{Interconnection between two sub-systems under perturbations.}
    \label{fig:smallgain}
    \vspace{-0.5cm}
\end{figure}
\begin{equation}\label{eq:sg_interconnection}
\begin{aligned}
x_1^+ &= g_1(x_1,x_2,u_1),\\
x_2^+ &= g_2(x_2,x_1,u_2),
\end{aligned}
\end{equation}
where $x_i\in\mathbb R^{n_i}$ and $u_i\in\mathbb R^{m_i}$, and $g_i(0,0,0)=0$ for $i\in\{1,2\}$.  Figure \ref{fig:smallgain} shows a block diagram representation of this system. 
\subsection{Main Assumptions and Result}
We make the following assumption on each subsystem:
\begin{assumption}\label{ass:sg}
Suppose that, for each $i\in\{1,2\}$, with $j=3-i$, there
exists a continuous function
$V_i:\mathbb{R}^{n_i}\to\mathbb{R}_{\geq0}$ satisfying the
following conditions:

\begin{enumerate}[(a)]
\item There exist
$\underline{\alpha}_i,\overline{\alpha}_i
\in\mathcal{K}_\infty$ such that
\begin{equation}
    \underline{\alpha}_i(|x_i|)
    \leq V_i(x_i)
    \leq \overline{\alpha}_i(|x_i|).
    \label{eq:sg_bounds}
\end{equation}

\item There exist $\rho_i\in[0,1)$ and
$\sigma_{ij},\sigma_i\in\mathcal{K}_\infty$ such that
\begin{equation}
    V_i(x_i^+)-V_i(x_i)
    \leq
    -(1-\rho_i) V_i(x_i)
    +\sigma_{ij}(V_j(x_j))
    +\sigma_i(|u_i|).
    \label{eq:sg_one_step}
\end{equation}

\item There exist $M\in\mathbb{Z}_{\geq1}$,
a continuous nondecreasing function
$\gamma_{ij}:\mathbb{R}_{\geq0}\to\mathbb{R}_{\geq0}$
satisfying $\gamma_{ij}(0)=0$, and
$\gamma_i\in\mathcal{K}_\infty$ such that every solution
satisfies
\begin{equation}
\begin{aligned}
    V_i(x_i(k+M))
    \leq {}&
    \gamma_{ij}\left(
        \max_{k\leq\ell<k+M}V_j(x_j(\ell))
    \right)\\
    &+
    \gamma_i\left(
        \|u_i\|_{[k,k+M-1],\infty}
    \right)
\end{aligned}
\label{eq:sg_finite_horizon}
\end{equation}
for every $k\in\mathbb{Z}_{\geq0}$. \hfill $\QEDB$
\end{enumerate}
\end{assumption}

\vspace{0.1cm}
To state a FxT-ISS result for the interconnected system \eqref{eq:sg_interconnection} based on finite-horizon FxT-ISS Lyapunov functions, we define the following loop gain maps for each $r\geq 0$:
\begin{align}
    \Gamma_{1,r}(s)
    &:=
    \gamma_{12}\!\left(
        \gamma_{21}(s)+\gamma_2(r)
    \right)+\gamma_1(r),
    \label{eq:Gamma1r}\\
    \Gamma_{2,r}(s)
    &:=
    \gamma_{21}\!\left(
        \gamma_{12}(s)+\gamma_1(r)
    \right)+\gamma_2(r).
    \label{eq:Gamma2r}
\end{align}
%
%
%\begin{assumption}\label{ass:fh_small_gain}
%There exist $N\in\mathbb{Z}_{\geq 1}$ and
%$\vartheta_1,\vartheta_2\in\mathcal{K}_\infty$ such that $\Gamma_{i,r}^{N}(s)\leq \vartheta_i(r)$ for all $(s,r)\in\mathbb{R}_{\geq0}^2,\quad i\in\{1,2\}.$
%
%\begin{equation}
%    \Gamma_{i,r}^{N}(s)
%    \leq \vartheta_i(r),
%    \qquad
%    \forall (s,r)\in\mathbb{R}_{\geq0}^2,
%    \quad i\in\{1,2\}.
%    \label{eq:fh_small_gain}
%\end{equation}
%\end{assumption}
%
The following theorem is the second main result of this paper.
\begin{theorem} \label{thm:fh_small_gain}
Suppose that Assumption \ref{ass:sg} holds. Then, the origin of~\eqref{eq:sg_interconnection} is FxT-ISS if there exist $N\in\mathbb{Z}_{\geq 1}$ and
$\vartheta_1,\vartheta_2\in\mathcal{K}_\infty$ such that 
\begin{equation}\label{smallgain2}
\Gamma_{i,r}^{N}(s)\leq \vartheta_i(r),~~\forall i\in\{1,2\}
\end{equation}
for all $(s,r)\in\mathbb{R}_{\geq0}^2$. \QEDB
\end{theorem}
\subsection{Application: Feedback Optimization with
Finite-Time Plant Settling}
\label{sec:feedback_optimization}

Consider the linear plant
\begin{equation}
x^+=Ax+Bv+d,
\label{eq:fo_plant}
\end{equation}
where $x\in\mathbb R^n$ is the state, $v\in\mathbb R^m$ is the control input, and and $d\in\mathbb R^n$ is an external disturbance. Suppose $A^{M_p}=0$ for some $M_p\in\mathbb Z_{\geq1}$. For a constant input $v$ and zero disturbance,
the equilibrium state and output are $x_{\rm eq}(v)=Xv$ and $y_{\rm eq}(v)=Sv$, where $X:=(I-A)^{-1}B$ and $S:=CX$. We consider the steady-state optimization problem
\begin{equation}
\begin{aligned}
\min_{v\in\mathbb R^m}\quad
\phi(v)
&:=\frac12v^\top Rv+b^\top v
  +\frac12(Sv)^\top Q_y(Sv)\\
&=\frac12v^\top Hv+b^\top v,
\end{aligned}
\label{eq:steady_state_optimization}
\end{equation}
where $R=R^\top\succ0$, $Q_y=Q_y^\top\succeq0$,
and $H:=R+S^\top Q_yS\succ0$. The unique optimizer to \eqref{eq:steady_state_optimization} is $v^\star=-H^{-1}b$.

Let $q$ denote a lifted iteration index. During each interval
$j\in\{qM_p,\ldots,(q+1)M_p-1\}$, we hold $v(j)=v_q$ constant. At the end of the interval, we measure

\begin{equation}
y_q=Cx((q+1)M_p)+\nu_q,
\label{eq:fo_measurement}
\end{equation}
where $\nu_q$ is measurement noise. The control update is:
\begin{equation}
v_{q+1}
=v_q-H^{-1}\bigl(Rv_q+b+S^\top Q_y y_q\bigr).
\label{eq:fo_newton}
\end{equation}
Thus, the optimizer uses the measured plant output to evaluate the output-dependent term in the gradient.  To study this system, let $\xi_q:=x(qM_p)$. Since $A^{M_p}=0$, we have
\begin{equation}
\xi_{q+1}=Xv_q+D_q,
\label{eq:fo_lifted_plant}
\end{equation}
where $D_q:=\sum_{\ell=0}^{M_p-1}$, $A^{M_p-1-\ell}d(qM_p+\ell)$. Consequently, $y_q=Sv_q+CD_q+\nu_q$, and $Rv_q+b+S^\top Q_y y_q
=\nabla\phi(v_q)+w_q$, with
\begin{equation}
w_q:=S^\top Q_y(CD_q+\nu_q).
\label{eq:fo_gradient_error}
\end{equation}
The settling interval therefore removes the initial
plant-state contribution from the measured gradient.

Let
$e_o(q):=v_q-v^\star$ and
$e_p(q):=\xi_q-Xv_q$.
Using \eqref{eq:fo_newton}--\eqref{eq:fo_gradient_error},
the lifted error dynamics are given by:
\begin{align}
e_o(q+1)&=-H^{-1}w_q,
\label{eq:fo_interconnection_o}\\
e_p(q+1)&=Xe_o(q)+XH^{-1}w_q+D_q.
\label{eq:fo_interconnection_p}
\end{align}
Note that although the implementation uses output feedback,
the lifted error dynamics have a cascade structure. To use Theorem \ref{thm:fh_small_gain}, let $V_o(e_o):=\|e_o\|$ and $V_p(e_p):=\|e_p\|$.
Using the common external input $u=(w,D)$ for both
subsystems, we obtain:
\begin{align*}
V_o(e_o^+)&\leq\|H^{-1}\|\|u\|,\\
V_p(e_p^+)&\leq\|X\|V_o(e_o)
 +(\|XH^{-1}\|+1)\|u\|.
\end{align*}
Thus, condition \eqref{eq:sg_one_step} holds with
$\rho_o=\rho_p=0$ and
\[
\begin{aligned}
\sigma_{op}(s)&=s,
&\sigma_o(s)&=\|H^{-1}\|s,\\
\sigma_{po}(s)&=(1+\|X\|)s,
&\sigma_p(s)&=(\|XH^{-1}\|+1)s.
\end{aligned}
\]
All these gains belong to $\mathcal K_\infty$. Condition \eqref{eq:sg_finite_horizon} holds with $M=1$ and
\[
\begin{aligned}
\gamma_{op}(s)&=0,
&\gamma_o(s)&=\|H^{-1}\|s,\\
\gamma_{po}(s)&=\|X\|s,
&\gamma_p(s)&=(\|XH^{-1}\|+1)s.
\end{aligned}
\]
The corresponding loop maps satisfy
\[
\begin{aligned}
\Gamma_{o,r}(s)&=\gamma_o(r),\\
\Gamma_{p,r}(s)&=\gamma_{po}(\gamma_o(r))+\gamma_p(r).
\end{aligned}
\]
Both expressions are independent of $s$.
Hence, the small-gain condition holds with $N=1$,
$\vartheta_o(r)=\|H^{-1}\|r$, and $\vartheta_p(r)
=\bigl(\|X\|\|H^{-1}\|+\|XH^{-1}\|+1\bigr)r$. Theorem~\ref{thm:fh_small_gain} therefore establishes
FxT-ISS of the lifted error system with horizon
$2NM=2$.

%To relate this result to the physical disturbances,
%suppose $\|d(k)\|\leq\bar d$ and
%$\|\nu_q\|\leq\bar\nu$. Then
%\[
%\begin{aligned}
%\|D_q\|&\leq\bar D
%:=\left(\sum_{\ell=0}^{M_p-1}\|A^\ell\|\right)\bar d,\\
%\|w_q\|&\leq\bar w
%:=\|S^\top Q_yC\|\bar D+\|S^\top Q_y\|\bar\nu.
%\end{aligned}
%\]
%These bounds also show that the lifted FxT-ISS
%estimate can be expressed in terms of the original
%plant and measurement disturbances.
%
%For explicit residual estimates,
%\eqref{eq:fo_interconnection_o} gives
%\[
%\|e_o(q)\|\leq\|H^{-1}\|\bar w,\qquad q\geq1.
%\]
%Moreover, for $q\geq1$,
%\[
%e_p(q+1)
%=XH^{-1}(w_q-w_{q-1})+D_q,
%\]
%and therefore
%\begin{equation}
%\|e_p(q)\|
%\leq2\|XH^{-1}\|\bar w+\bar D,
%\qquad q\geq2.
%\label{eq:fo_residual}
%\end{equation}
%In the absence of disturbances, $e_o(q)=0$ for
%$q\geq1$ and $e_p(q)=0$ for $q\geq2$.
%Thus, the optimizer reaches $v^\star$ after $M_p$
%physical plant steps, and the plant reaches
%$Xv^\star$ after at most $2M_p$ steps.

Figure 2 shows a numerical example of the proposed controller, where
\[
A=\begin{bmatrix}0&1&0\\0&0&1\\0&0&0\end{bmatrix},
\quad
B=\begin{bmatrix}0\\0\\1\end{bmatrix},
\quad
C=\begin{bmatrix}1&0&0\end{bmatrix},
\]
and $M_p=3$, $R=Q_y=1$, and $b=-2$. In this case, we have $X=(1,1,1)^\top$, $S=1$, $H=2$, and $v^\star=1$. We used $x(0)=(3,-2,4)^\top$ and $v_0=4$. For disturbance bounds, we used $\bar d=0.03$ and $\bar\nu=0.02$. The estimates above give
$\bar D=0.09$, $\bar w=0.11$, and
\[
\|e_o(q)\|\leq0.055\quad(q\geq1),~~
\|e_p(q)\|\leq0.09+0.11\sqrt3\quad(q\geq2).
\]
Note that the nominal and disturbed trajectories are compared
using the lifted iteration index $q$.
\begin{figure}[t]
    \centering
    \includegraphics[width=0.9\columnwidth]{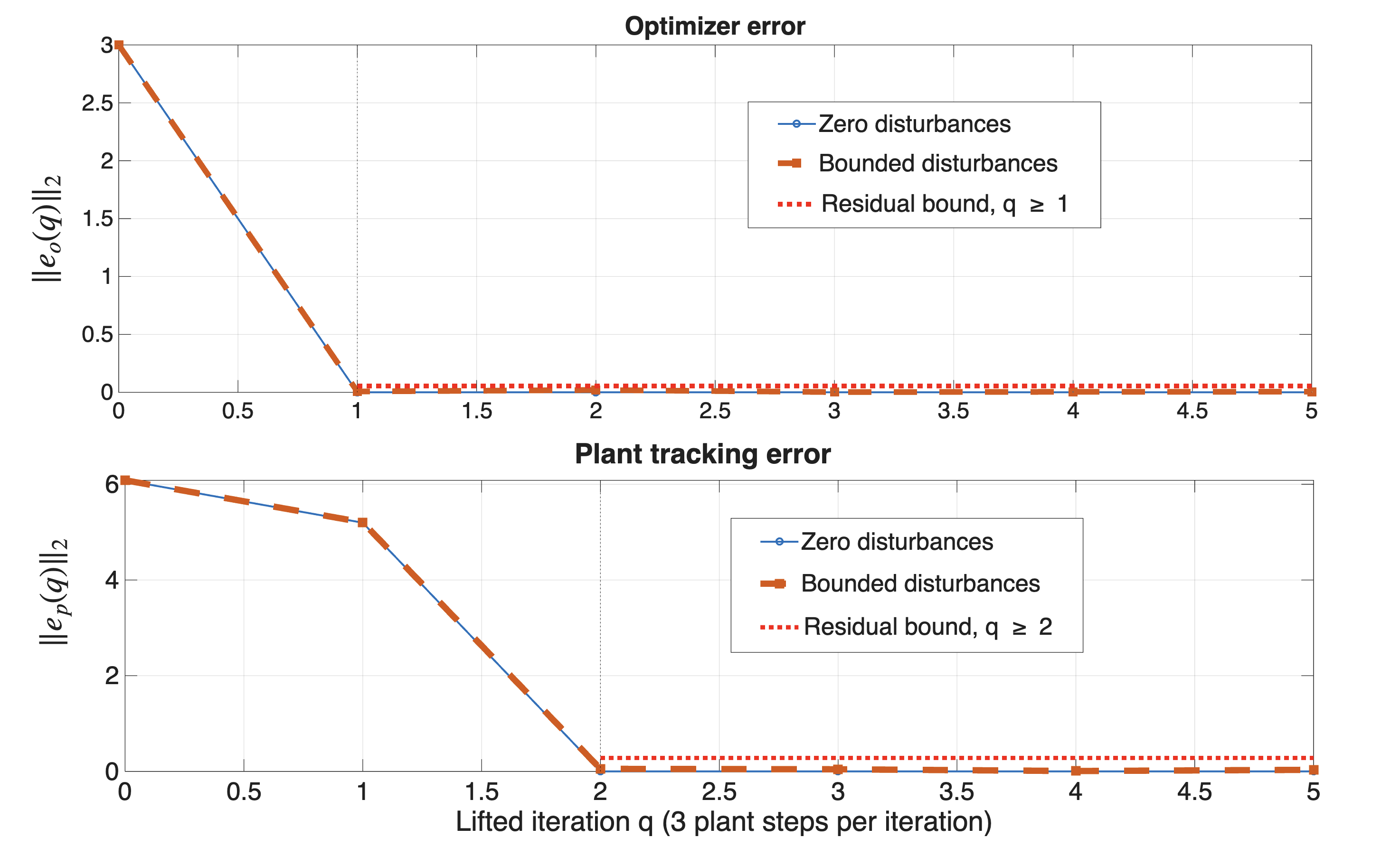}
    \caption{Plant with finite settling time driven by Newton-based control.}
    \label{fig:example2}
    \vspace{-0.5cm}
\end{figure}
\section{Proofs}
\label{secproofs}
\subsection{Proof of Proposition \ref{equivalentprop}:} It suffices to establish the equivalence between conditions \eqref{eq:fh_dissipation} and \eqref{implicationcondition}. To show necessity, suppose that $V$ is a finite-horizon FxT-ISS Lyapunov function. Then, by definition, there exist $\rho_0\in[0,1)$ and $\sigma\in\mathcal K_\infty$ such that \eqref{eq:fh_dissipation} holds. Choose any $\rho\in(\rho_0,1)$ and define $\chi(r):=
\frac{\sigma(r)}{\rho-\rho_0}\in\mathcal K_\infty$. If $V(x)\geq\chi(\|u\|)$, then $\sigma(\|u\|)\leq(\rho-\rho_0)V(x)$. Consequently,
\begin{equation}
V(x^+)\leq\rho_0V(x)+(\rho-\rho_0)V(x)=\rho V(x).
\end{equation}
Thus, the implication-type condition holds.

To show sufficiency, suppose now that conditions \eqref{eq:fh_comparison} and \eqref{finitehorizonbound} hold. If $V(x)\geq\chi(\|u\|)$, then condition \eqref{implicationcondition} directly gives $V(x^+)\leq\rho V(x)$. On the other hand, if $V(x)<\chi(\|u\|)$, define for $r\geq0$: 
\[
\widehat{\sigma}(r)
:=
\max\left\{
V(g(x,u)):
V(x)\leq\chi(r),\ \|u\|\leq r
\right\}.
\]
The function $\widehat{\sigma}(r)$ is well defined and finite. Indeed, from \eqref{eq:fh_comparison} we have that $V(x)\leq\chi(r)$ implies $|x|_{\mathcal A}\leq\alpha_1^{-1}(\chi(r))$. Since $\mathcal A$ is compact, the set $\left\{x\in\mathbb R^n:|x|_{\mathcal A}\leq\alpha_1^{-1}(\chi(r))
\right\}$ is compact. Hence, the set $\left\{
(x,u):
V(x)\leq\chi(r),\ \|u\|\leq r
\right\}$ is compact. Since $g$ and $V$ are continuous, the maximum defining
$\widehat{\sigma}(r)$ exists and is finite. Moreover, $\widehat{\sigma}(0)=0$. Indeed, $\chi(0)=0$, and therefore $V(x)\leq\chi(0)$ implies $V(x)=0$. Applying condition \eqref{implicationcondition} with $u=0$ gives
\[
V(g(x,0))
\leq
\rho V(x)=0.
\]
Thus, $V(g(x,0))=0$. By continuity of $g$ and $V$, compactness of $\mathcal A$, and the comparison bounds, it follows that $\widehat{\sigma}(r)\to0
\qquad\text{as }r\to0^+$. Furthermore, $\widehat{\sigma}$ is nondecreasing and finite on
$\mathbb R_{\geq0}$. Hence, there exists
$\sigma\in\mathcal K_\infty$ such that $\widehat{\sigma}(r)\leq\sigma(r)$, for all $r\geq0$. Now, if $V(x)<\chi(\|u\|)$, then, taking $r=\|u\|$, we get
\[
V(g(x,u))
\leq
\widehat{\sigma}(\|u\|)
\leq
\sigma(\|u\|).
\]
Combining the two cases gives $V(g(x,u))\leq\max\left\{\rho V(x),\sigma(\|u\|)\right\}$. Finally, using $\max\{a,b\}\leq a+b$ for all $a,b\geq0$, we obtain $V(g(x,u))\leq\rho V(x)+\sigma(\|u\|)$. Thus, $V$ is a finite-horizon FxT-ISS Lyapunov function. \hfill $\blacksquare$

\subsection{Proof of Theorem \ref{theoremLyapu} :}  Let $V_j:=V(x(j))$ and $U_j:=\|u\|_{[0,j-1],\infty}$. From \eqref{eq:fh_dissipation}, we have
$V_{j+1}\leq\rho V_j+\sigma(\|u(j)\|)$. Iterating this inequality yields $V_j\leq \rho^jV_0+\sum_{k=0}^{j-1}\rho^{j-1-k}\sigma(\|u(k)\|)$.
Since $\sigma$ is increasing,
\begin{equation}
V_j\leq \rho^j V_0+\left(\sum_{k=0}^{j-1}\rho^k\right)\sigma(U_j).
\end{equation}
For $0\leq j<M$, define $C_M:=\sum_{k=0}^{M-1}\rho^k=\frac{1-\rho^M}{1-\rho}$. Then $V_j\leq\rho^j V_0+C_M\sigma(U_j),~~0\leq j<M$. Using \eqref{eq:fh_comparison}, we have $V_0\leq\alpha_2(|x(0)|_{\mathcal A})$ and hence
\begin{equation}
\label{eq:pre_horizon_bound}
V_j
\leq
\rho^j\alpha_2(|x(0)|_{\mathcal A})
+
C_M\sigma(U_j),
\qquad 0\leq j<M.
\end{equation}

For $j\geq M$, applying \eqref{finitehorizonbound} at time $j-M$ gives
\begin{equation}
V_j=V\!\left(G_M(x(j-M),u(j-M),\ldots,u(j-1))\right)
\end{equation}
and therefore
\begin{equation}
V_j\leq\gamma\!\left(
\max_{j-M\leq k\leq j-1}\|u(k)\|
\right).
\end{equation}
Since $\max_{j-M\leq k\leq j-1}\|u(k)\|\leq U_j$, and $\gamma$ is increasing, we get $V_j\leq\gamma(U_j)$ for all $j\geq M$.
%
%\begin{equation}
%\label{eq:post_horizon_bound}
%V_k\leq\gamma(U_k),
%\qquad k\geq M.
%\end{equation}
%
Choose $\bar\rho\in(\rho,1)$ and define
\[
b(s,j):=
\begin{cases}
\bar\rho^j\alpha_2(s), & 0\leq j<M,\\
0, & j\geq M.
\end{cases}
\]
and $\eta(r):=\max\left\{C_M\sigma(r),\gamma(r)\right\}$. Combining \eqref{eq:pre_horizon_bound} and $V_j\leq\gamma(U_j)$ for all $j\geq M$ gives, 
\begin{equation}
V_j\leq b(|x(0)|_{\mathcal A},j)+\eta(U_j),~~\forall~j\geq0.
\end{equation}
Using the lower bound in \eqref{eq:fh_comparison}, we get $\alpha_1(|x(j)|_{\mathcal A}) \leq b(|x(0)|_{\mathcal A},j)+\eta(U_j)$. Therefore,
\[
|x(j)|_{\mathcal A}
\leq
\alpha_1^{-1}\!\left(
b(|x(0)|_{\mathcal A},j)+\eta(U_j)
\right).
\]
Using the inequality $\alpha_1^{-1}(a+b)
\leq
\max\left\{
\alpha_1^{-1}(2a),\alpha_1^{-1}(2b)
\right\}$, for $a,b\geq0$, we obtain
\begin{equation}
|x(j)|_{\mathcal A}
\leq
\max\left\{
\beta(|x(0)|_{\mathcal A},j),\chi(U_j)
\right\},
\end{equation}
where
$\beta(s,j):=\alpha_1^{-1}(2b(s,j))$
and $\chi(r):=\alpha_1^{-1}(2\eta(r))$.
Since $0<\bar\rho<1$, the function $\beta(\cdot,j)$
belongs to $\mathcal K$ for every $j<M$,
$\beta(s,\cdot)$ is nonincreasing, and
$\beta(s,j)=0$ for every $j\geq M$.
Thus, $\beta\in\mathcal{KL}_{\mathrm{FxT}}$ and
$\chi\in\mathcal K_\infty$, which proves the claim. \hfill $\blacksquare$

\subsection{Proof of Proposition \ref{thm:extinction_implies_fh}}
It suffices to show that $V$ satisfies inequalities \eqref{eq:fh_dissipation} and \eqref{finitehorizonbound}.

We first show that $V$ satisfies \eqref{eq:fh_dissipation}. Since
$\psi\in\mathcal{E}_{\mathrm{FxT}}^{M}$, by \cite[Lemma 1]{ACC26A} there exist $\varepsilon>0$ and $H<\infty$ such that: (a) $\psi(s)\geq s$, for all $0\leq s\leq\varepsilon$, and (b) $\psi(s)\geq s-H$ for all
$s\geq0$. Moreover, (c) for every $s\geq 2H$, we have that $ \frac{\psi(s)}{s}
    \geq 1-\frac{H}{s}
    \geq \frac{1}{2}$.
%\begin{equation}
%    \psi(s)\geq s,
%    \qquad 0\leq s\leq\varepsilon,
%    \label{eq:psi_near_origin}
%\end{equation}
%and
%\begin{equation}
%    \psi(s)\geq s-H,
%    \qquad s\geq0.
%    \label{eq:psi_global_bound}
%\end{equation}
%\begin{equation}
%    \frac{\psi(s)}{s}
%    \geq 1-\frac{H}{s}
%    \geq \frac{1}{2}.
%    \label{eq:psi_large_ratio}
%\end{equation}
%
Moreover, since $\psi$ is continuous and positive definite, for each $s>0$ the map
$s\mapsto\psi(s)/s$ is continuous and strictly positive on every
compact subset of $(0,\infty)$. Hence, whenever
$\varepsilon<2H$, the quantity $c_0:=
    \min_{s\in[\varepsilon,2H]}
    \frac{\psi(s)}{s}$i s well defined and satisfies $c_0>0$. Combining this fact with
properties (a) and (c), there exists a constant $c\in(0,1]$ such that $\psi(s)\geq cs$, for all  $s\geq0$. For example, if $\varepsilon<2H$, one may take $c:=\min\left\{\frac{1}{2},c_0\right\}$. If $\varepsilon\geq2H$, one may simply take $c=1/2$. Let $\varrho:=1-c\in[0,1)$. Then, for every $s\geq0$,
\[
    \min\{\psi(s),s\}\geq cs=(1-\varrho)s.
\]
Therefore, the extinction-based dissipation inequality gives
\begin{align}
    V(x^+)-V(x)
    &\leq
    -\min\{\psi(V(x)),V(x)\}+\rho(|u|)
    \nonumber\\
    &\leq
    -(1-\varrho)V(x)+\rho(|u|).
    \label{eq:one_step_fh_from_extinction}
\end{align}
Thus, $V$ satisfies \eqref{eq:fh_dissipation}. To prove \eqref{finitehorizonbound}, define $g_{\psi}(s):=[s-\psi(s)]_{+}$, $h_{\psi}(s):=\max_{0\leq \tau\leq s}g_{\psi}(\tau)$. Since $\psi\in\mathcal{E}_{\mathrm{FxT}}^{M}$, it follows that
$h_{\psi}^{M}(s)=0$ for all $s\geq0$. Fix
$k\in\mathbb{Z}_{\geq0}$ and let $v_j:=V(x(k+j))$, $r:=\|u\|_{[k,k+M-1],\infty}$. For $j\in\{0,\ldots,M-1\}$, the extinction-based
dissipation inequality gives
\begin{align}
v_{j+1}
&\leq
v_j-\min\{\psi(v_j),v_j\}+\rho(r)=g_{\psi}(v_j)+\rho(r)
\nonumber\\
&\leq
h_{\psi}(v_j)+\rho(r).
\label{eq:scalar_recursion_fh}
\end{align}
For each $r\geq0$, define $T_r(s):=h_{\psi}(s)+\rho(r)$. Since $h_{\psi}$ is nondecreasing, repeated application of
\eqref{eq:scalar_recursion_fh} yields $v_M\leq T_r^M(v_0)$. The extinction property implies that $h_{\psi}$ is globally
bounded. Let $H:=\sup_{s\geq0}h_{\psi}(s)<\infty$. Hence,
$T_r(s)\leq H+\rho(r)$, for all $s\geq0$. Consequently, $T_r^M(s)\leq T_r^{M-1}(H+\rho(r))$, for all $s\geq0$. Define $\Gamma(r):=T_r^{M-1}(H+\rho(r))$. The map $\Gamma$ is continuous because $h_{\psi}$ and $\rho$
are continuous. Moreover, $\Gamma(0)=h_{\psi}^{M-1}(H)$. Since $H$ belongs to the range closure of $h_{\psi}$ and
$h_{\psi}^{M}\equiv0$, continuity of $h_{\psi}^{M-1}$ gives $h_{\psi}^{M-1}(H)=0$. Thus, $\Gamma(0)=0$, and therefore there exists
$\sigma\in\mathcal{K}_{\infty}$ such that $\Gamma(r)\leq\sigma(r)$,
for all $r\geq0$. It follows that
\[
V(x(k+M))
=
v_M
\leq
\sigma\!\left(
\|u\|_{[k,k+M-1],\infty}
\right),
\]
which is precisely inequality \eqref{finitehorizonbound}.
\hfill $\blacksquare$

\subsection{Proof of Proposition \ref{propnecessarycondition}}
To prove item (a), suppose that $V$ is an extinction-based FxT-ISS Lyapunov function. Then there exist $\alpha_1,\alpha_2\in\mathcal K_\infty$ and $\psi\in\mathcal E_{\mathrm{FxT}}^N$ such that \eqref{eq:general_lyap_bounds} holds. Setting $u=0$ in \eqref{eq:additive_dissipation}, we get
\begin{equation}
V(g(x))\leq[V(x)-\psi(V(x))]_+\leq h_\psi(V(x)),
\end{equation}
where $h_\psi(s):=\max_{0\leq r\leq s}[r-\psi(r)]_+$.
The map $h_\psi$ is nondecreasing and satisfies $h_\psi^N\equiv0$. Therefore, it is bounded: otherwise, every finite composition of $h_\psi$ would be unbounded. Let $H:=\sup_{s\geq0}h_\psi(s)<\infty$. It follows that $\alpha_1(|g(x)|_{\mathcal A})\leq V(g(x))\leq H$ for all $x$. Hence, $|g(x)|_{\mathcal A}\leq\alpha_1^{-1}(H)$, which proves item (a).

To prove item (b), let $\rho\in(0,1)$ and define $V(x):=\sum_{\ell=0}^{M-1}\rho^{-\ell}\|g^\ell(x)\|$. Continuity of $g$ and $g(0)=0$ imply that $V$ is continuous and $V(0)=0$. Moreover, $V(x)\geq\|x\|$. Define
\[
b(s):=\max_{\|z\|\leq s}V(z),\qquad s\geq0.
\]
Continuity of $V$ implies that $b$ is finite, continuous, and nondecreasing, with $b(0)=0$. Consequently, $\alpha_2(s):=s+b(s)$ belongs to $\mathcal K_\infty$, and $\|x\|\leq V(x)\leq\alpha_2(\|x\|)$.
Using $g^M(x)\equiv0$, we obtain
\begin{align*}
V(g(x))
&=\sum_{\ell=0}^{M-1}\rho^{-\ell}\|g^{\ell+1}(x)\|=\rho\sum_{\ell=1}^{M}\rho^{-\ell}\|g^\ell(x)\|\\
&=\rho\bigl(V(x)-\|x\|\bigr)\leq\rho V(x).
\end{align*}
which establishes \eqref{eq:fh_dissipation} (with $u=0$).  Furthermore, the $M$-step transition map satisfies $G_M(x,u_0,\ldots,u_{M-1})=g^M(x)=0$. Therefore, $V(G_M(x,u_0,\ldots,u_{M-1}))=0
\leq 0$, which establishes \eqref{finitehorizonbound} (wtih $u=0$) \hfill $\blacksquare$

\subsection{Proof of Theorem \ref{thm:fh_small_gain}}
Let $v_i(k):=V_i(x_i(k))$, $i\in\{1,2\}$, and, for integers $0\leq k\leq L$, define $S_i^L(k):=\max_{k\leq\ell\leq L}v_i(\ell)$. Also, let $U_{k,L}:=\|u\|_{[k,L-1],\infty}$. Fix $L\geq k+M$. For every
$q\in\{k+M,\ldots,L\}$, applying
\eqref{eq:sg_finite_horizon} at time $q-M$ gives
\begin{align*}
v_i(q)
&\leq
\gamma_{ij}\left(
    \max_{q-M\leq\ell<q}v_j(\ell)
\right)
+\gamma_i\left(
    \|u_i\|_{[q-M,q-1],\infty}
\right)\\
&\leq
\gamma_{ij}(S_j^L(k))
+\gamma_i(U_{k,L}).
\end{align*}
Taking the maximum over
$q\in\{k+M,\ldots,L\}$ yields
\begin{equation}
    S_i^L(k+M)
    \leq
    \gamma_{ij}(S_j^L(k))
    +\gamma_i(U_{k,L}).
    \label{eq:tail_one_step}
\end{equation}
Suppose now that $L\geq k+2M$. Applying
\eqref{eq:tail_one_step} successively to both subsystems
gives
\begin{align*}
S_1^L(k+2M)
&\leq
\gamma_{12}(S_2^L(k+M))
+\gamma_1(U_{k,L})\\
&\leq
\gamma_{12}\left(
    \gamma_{21}(S_1^L(k))
    +\gamma_2(U_{k,L})
\right)
+\gamma_1(U_{k,L})\\
&=
\Gamma_{1,U_{k,L}}(S_1^L(k)).
\end{align*}
Similarly, $S_2^L(k+2M)\leq \Gamma_{2,U_{k,L}}(S_2^L(k))$. Iterating these inequalities $N$ times yields, whenever
$L\geq k+2NM$,
\begin{equation}
    S_i^L(k+2NM)
    \leq
    \Gamma_{i,U_{k,L}}^N(S_i^L(k))
    \leq
    \vartheta_i(U_{k,L}),
    \qquad i\in\{1,2\},
    \label{eq:tail_N}
\end{equation}
where the last inequality follows from \eqref{smallgain2}. Let $K:=2NM$. For any $j\geq K$, choose $k=0$ and $L=j$ in \eqref{eq:tail_N}. Since $j\in[K,j]$, it follows that
\begin{equation}
    v_i(j)
    \leq
    S_i^j(K)
    \leq
    \vartheta_i(U_j),
    \qquad
    U_j:=\|u\|_{[0,j-1],\infty}.
    \label{eq:post_horizon_bound}
\end{equation}
Thus, after the uniform horizon $K$, the bound is
independent of the initial condition. In particular, for
$u=0$, since $\vartheta_i(0)=0$, $v_i(j)=0$, for $j\geq K$. It remains to bound the transient $0\leq j<K$. Define $R(j):=\max\{v_1(j),v_2(j)\}$, and let $\rho:=\max\{\rho_1,\rho_2\}$, $\sigma(s):=\max\{\sigma_{12}(s),\sigma_{21}(s)\}$, $\sigma_u(s):=
    \max\{\sigma_1(s),\sigma_2(s)\}$. It follows from \eqref{eq:sg_one_step} that
\begin{equation}
    R(j+1)
    \leq
    \rho R(j)+\sigma(R(j))+\sigma_u(|u(j)|).
    \label{eq:R_transient}
\end{equation}
Moreover, define $\overline{\alpha}(s):=
\max\left\{
\overline{\alpha}_1(s),
\overline{\alpha}_2(s)
\right\}$. Then, condition~\eqref{eq:sg_bounds} gives
$R(0)\leq \overline{\alpha}(|x(0)|)$. Since $K$ is finite, repeated application of
\eqref{eq:R_transient} yields, for each $0\leq j<K$,
functions $p_j\in\mathcal{K}_\infty$ and
$q_j\in\mathcal{K}$ such that
\begin{equation}
    R(j)
    \leq
    p_j(|x(0)|)+q_j(U_j).
    \label{eq:R_finite_bound}
\end{equation}
Define $\underline{\alpha}(s):=
    \min_{i\in\{1,2\}}
    \underline{\alpha}_i
    \left(\frac{s}{\sqrt{2}}\right)$. Since at least one component satisfies
$|x_i(j)|\geq |x(j)|/\sqrt{2}$, we have $\underline{\alpha}(|x(j)|)\leq R(j)$. Hence, using \eqref{eq:R_finite_bound} and $\underline{\alpha}^{-1}(a+b) \leq \underline{\alpha}^{-1}(2a)+ \underline{\alpha}^{-1}(2b)$, we obtain
\begin{equation}
    |x(j)|
    \leq
    b_j(|x(0)|)+c_j(U_j),
    \qquad 0\leq j<K,
    \label{eq:state_transient}
\end{equation}
where $b_j(s):= \underline{\alpha}^{-1}(2p_j(s))$, $c_j(s):=
    \underline{\alpha}^{-1}(2q_j(s))$. For $0\leq j<K$, define $\beta(s,j):=
    \max_{j\leq\ell<K}b_\ell(s)$, and set $\beta(s,j):=0$, for  $j\geq K$. By construction, $\beta(\cdot,j)\in\mathcal{K}$ for
$j<K$, $\beta(s,\cdot)$ is nonincreasing, and
$\beta(s,j)=0$ for every $s\geq0$ and $j\geq K$.
Therefore, $\beta\in\mathcal{KL}_{\rm FxT}$. To construct the input gain, define $c(s):=\max_{0\leq j<K}c_j(s).$ Moreover, from \eqref{eq:post_horizon_bound}, we have $|x_i(j)|
    \leq
    \underline{\alpha}_i^{-1}
    \bigl(\vartheta_i(U_j)\bigr)$, for $j\geq K$. Consequently, $|x(j)|
    \leq d(U_j)$ for $j\geq K$, where $d(s):=
    \left(
        \sum_{i=1}^2
        \left[
            \underline{\alpha}_i^{-1}
            \bigl(\vartheta_i(s)\bigr)
        \right]^2
    \right)^{1/2}$. Let $\chi\in\mathcal{K}$ be any class-$\mathcal{K}$
majorant of $\max\{c,d\}$. Combining
\eqref{eq:state_transient} with the post-horizon estimate
then gives
\[
    |x(j)|
    \leq
    \beta(|x(0)|,j)+\chi(U_j),
    \qquad j\in\mathbb{Z}_{\geq0}.
\]
Finally, using $a+b\leq\max\{2a,2b\}$ for  $a,b\geq0$, and defining $\widetilde{\beta}(s,j):=2\beta(s,j)$, and $\widetilde{\chi}(s):=2\chi(s)$,
we obtain
\begin{equation}
    |x(j)|
    \leq
    \max\left\{
        \widetilde{\beta}(|x(0)|,j),
        \widetilde{\chi}
        \left(\|u\|_{[0,j-1],\infty}\right)
    \right\}.
\end{equation}
Since $\widetilde{\beta}\in\mathcal{KL}_{\rm FxT}$, the origin of~\eqref{eq:sg_interconnection} is FxT-ISS. \hfill $\blacksquare$

\section{Conclusions}
\label{sec:conclusions}
This paper introduced a finite-horizon Lyapunov framework for studying fixed-time ISS in discrete-time systems. The proposed conditions capture the uniform loss of dependence on the initial state after a finite number of iterations, while allowing this property to arise from the multi-step structure of the dynamics. We also established a finite-horizon small-gain theorem for interconnected systems, providing conditions under which this property is preserved under feedback interconnections. The framework was illustrated through examples involving deadbeat dynamics, distributed consensus, and finite-step optimization. These tools complement the extinction-based Lyapunov approach studied in the companion paper \cite{ACC26A}, and open the door for a systematic analysis and design of FxT-ISS algorithms in discrete-time.

\bibliographystyle{IEEEtran}
\bibliography{Biblio.bib}

\end{document}